\documentclass[11pt]{amsart}

\newtheorem{theorem}{Theorem}

\newtheorem{lemma}{Lemma}
\newtheorem{corollary}{Corollary}

\newenvironment{pf*}[1]{\medskip \noindent {\em #1.} }{\endproof \medskip}

\newcommand{\zz}[1]{\mathbb #1}

\title[Critical Scaling for the SIS Epidemic]{Critical Scaling for the
Simple SIS Stochastic Epidemic} \author{R. G. Dolgoarshinnykh}
\author{Steven P. Lalley}

\thanks{S.L. supported by NSF grant DMS-0405102}
\date{\today}

\begin{document}

\maketitle

\begin{abstract}
We exhibit a scaling law for the critical SIS stochastic epidemic: If
at time $0$ the population consists of $\sqrt{N}$ infected and
$N-\sqrt{N}$ susceptible individuals, then when time and number
currently infected are both scaled by $\sqrt{N}$, the resulting
process converges, as $N \rightarrow \infty$, to a diffusion process
related to the Feller diffusion by a change of drift. As a
consequence, the rescaled size of the epidemic has a limit law that
coincides with that of a first-passage time for the standard
Ornstein-Uhlenbeck process. These results are the analogues for the
SIS epidemic of results of Martin-L\"{o}f \cite{mlof} for the simple
SIR epidemic.
\end{abstract}

\section{Introduction}
\label{sec:intro}

Among the most thoroughly studied stochastic epidemic models are the
simple SIR, and SIS epidemics (see Weiss and Dishon
\cite{weiss-dishon} for the origin of the SIS model); and among the
many problems associated with these models perhaps the most basic and
most interesting have to do with the duration and size of the
epidemic. When the epidemic is either sub- or supercritical, the
large-population behavior of the duration and size is reasonably
well-understood: in the subcritical case, the epidemic is
stochastically dominated by a subcritical Galton-Watson process (see
below), and in the supercritical case the epidemic may become
\emph{endemic} (see Norden \cite{norden} and Kryscio and Lef\`{e}vre
\cite{k-l}). For critical epidemics, large-population asymptotics are
more delicate. In 1998, Martin-L\"{o}f \cite{mlof} (see also Aldous
\cite{aldous}) showed that the size $S=S_{N}$ has an interesting and
nontrivial asymptotic behavior as the population size $N \rightarrow
\infty$: If the number $X_{0}$ of individuals initially infected is of
order $bN^{1/3}$, then the size $S_{N}$ (defined to be the sum of the
total infection times over the whole population) has a limit
distribution
\begin{equation}\label{eq:ml}
    S_{N}/N^{2/3} \stackrel{\mathcal{D}}{\longrightarrow} T^{*}_{b}
\end{equation}
where $T^{*}_{b}$ is the first passage time of $W (t)+t^{2}/2$ to the
level $b$, and $W(t)$ is a standard Wiener process. Furthermore, $1/3$
is the critical exponent, in that the quadratic drift is not felt if
$X_{0}$ is much smaller than $N^{1/3}$: If $X (0)\sim bN^{\alpha}$ for
some $\alpha<1/3$ then $S_{N}/N^{2\alpha}$ converges in law to the
first passage time of $W(t)$ to the level $b$. (Martin-L\"{o}f neither
proved nor stated this, but it can be deduced from his methods.)

The purpose of this note is to establish
an analogous scaling law for the critical SIS epidemic: If the number
$X_{0}$ of individuals initially infected is of order $b\sqrt{N}$,
then as $N \rightarrow \infty$,
\begin{equation}\label{eq:scalingLaw}
    S_{N}/N \stackrel{\mathcal{D}}{\longrightarrow} \tau_{b},
\end{equation}
where $\tau_{b}$ is the time of first passage to $0$ by a standard
Ornstein-Uhlenbeck process started at $b$. Our approach has a rather
different character than those of Aldous and Martin-L\"{o}f. We shall
establish that the SIS epidemic process itself, suitably rescaled,
converges in law to a diffusion process we dub the \emph{attenuated
Feller diffusion}. This has a law absolutely continuous with respect to,
but not equal to, that of the standard {Feller diffusion}. (See
equation (\ref{eq:feller}) below for the stochastic differential
equation governing the Feller diffusion, and
(\ref{eq:attenuatedFellerDiffusion}) for the attenuated Feller
diffusion.) Furthermore, we will show that $1/2$ is the critical
exponent, in the following sense: If $X (0)\sim bN^{\alpha}$ for some
$\alpha <1/2$ then the rescaled SIS process converges in law to a
standard Feller diffusion with \emph{no} drift. It will follow that
the duration of the epidemic, rescaled by $\sqrt{X (0)}$, converges in
law as $N \rightarrow \infty$ to the first passage time to $0$ of the
corresponding Feller diffusion.

Our analysis will show that the critical \emph{scaling window} for the
transmissivity parameter $\beta$ (see below) is on the order
$1/\sqrt{N}$. This scaling window has also been observed --- but in a
different context --- by Nasell \cite{nasellA}, \cite{nasellB}. Nasell
shows that the quasistationary distribution of the SIS epidemic
undergoes a scaling transition when the  transmissivity parameter
varies from below $1-O (N^{-1/2})$ to above $1+O (N^{-1/2})$. This
phenomenon does not seem to be directly linked to the critical scaling
in our Theorem \ref{theorem:main}.

\section{The SIS Epidemic and its Branching Envelope}\label{sec:sis}

\subsection{SIS Model}\label{ssec:sis} The SIS epidemic  is a
continuous-time birth-death Markov chain $X_{t}=X (t)$ on the state
space $[N]:=\{0,1,2,\dots ,N \}$ whose infinitesimal transition
probabilities are as follows:
\begin{align}\label{eq:rates}
    P\{X (t +\delta t)&=x+1 \, | \, X (t)=x \}
    = \beta x (1-x/N) \delta t +o (x\delta t),\\
\notag    P\{X (t +\delta t)&=x-1 \, | \, X (t)=x \}
    = x\delta t  + o (x\delta t),\\
\notag    P\{X (t +\delta t)&=x \, | \, X (t)=x \}
    = 1- \beta x (1-x/N) \delta t -x\delta t +o (x\delta t)).
\end{align}
These describe a population of $N$ individuals in which $X_{t}$ are
infected and the remainder $N-X_{t}$ are susceptible to infection at
time $t$. Infected individuals recover at rate $1$, after which they
once again become susceptible to (re-)infection, and susceptible
individuals become infected at rate $\beta X_{t}/N$ proportional to
the number of infected individuals in the population. The epidemic
ends at the first time $T=T_{N}=t$ when $X_{t}=0$ (note that state $0$
is absorbing). The epidemic is said to be \emph{critical} when $\beta
=1$, and \emph{nearly critical} when $\beta =1+\lambda /\sqrt{N}$.

\subsection{Branching Envelope}\label{ssec:gwe} When the number of
individuals infected is small compared to the population size, the
epidemic evolves approximately as a continuous-time branching process
$Z (t)=Z_{t}$ with infinitesimal transition probabilities
\begin{align}\label{eq:gwe}
    P\{Z (t+\delta t)&=z+1 \, | \, Z (t)=z\}
    = \beta z \delta t +o (z \delta  t),\\
\notag  P\{Z (t+\delta t)&=z-1 \, | \, Z (t)=z\}
    = z \delta t +o (z \delta  t),\\
\notag  P\{Z (t+\delta t)&=z \, | \, Z (t)=z\}
    = 1 - (\beta +1)z \delta  t +o (z \delta  t).
\end{align}
We shall refer to this process as the \emph{branching envelope} of the
SIS process. Observe that the death rate $x$ is the same as for the
SIS epidemic, but the the birth rate $\beta x$ dominates the birth
rate $\beta x (1-x/N)$ of the SIS process; the difference $\beta
x^{2}/N$ will be called the \emph{attenuation} or \emph{attrition
rate}.  It is possible, by a standard construction, to build the SIS
process $X (t)$ and its branching envelope $Z (t)$ on the same
probability space in such a way that $X (0)=Z (0)$ and $X (t)\leq Z
(t)$ for all $t\geq 0$. Thus, the size and duration of the SIS
epidemic are stochastically dominated by the total progeny and
extinction time of the branching envelope.

\subsection{Critical Scaling for the Branching
Envelope}\label{ssec:fellerLimit} It was proved by
Feller~\cite{feller} that a critical branching process, when
properly renormalized, behaves approximately as a \emph{Feller
diffusion} with drift $\lambda Y_{t}$, that is, a solution of the
stochastic differential equation
\begin{equation}\label{eq:feller}
    dY_{t}= \lambda Y_{t} \,dt +\sqrt{Y_{t}}dW_{t}
\end{equation}
where $W_{t}$ is a standard Wiener process. (Equivalently, the Feller
diffusion with drift parameter $\lambda$ may be described as the
diffusion process on $[0,\infty )$ with infinitesimal generator
$\mathcal{G}^{\lambda}=\lambda x \partial_{x}+x\partial^{2}_{xx}/2$.)
Feller's theorem (see Jirina \cite{jirina} and Lindvall
\cite{lindvall} for the proof) asserts that if $Z^{m}(t)$ is a
sequence of branching processes satisfying (\ref{eq:gwe}) with $\beta
=\beta_{m}=1+\lambda /m$ and with $Z^{m} (0)\sim bm$ for some $b>0$
and $\lambda \in \zz{R}$ then
\begin{equation}\label{eq:bpConvergence}
    Z^{m} (mt)/m  \stackrel{\mathcal{D}}{\longrightarrow} Y_{t}
\end{equation}
where $Y_{t}$ is the Feller diffusion with drift parameter $\lambda$
and initial value $Y_{0}=b$.

\subsection{Critical Scaling for the SIS
Process}\label{ssec:mainResult} Because the branching envelope
stochastically dominates the SIS process, the scaling law
(\ref{eq:bpConvergence}) limits the duration and growth of the
critical and near-critical SIS
process. Since time is scaled by the factor $m$, where $Z^{m} (0)\sim
bm$, it follows that the corresponding SIS started with $X (0)\sim
bN^{\alpha}$ infected individuals cannot have duration longer than
$O_{P} (N^{\alpha})$ time units. Consequently, we should expect that
if the attenuation rate, divided by the scale factor $N^{\alpha}$ and
integrated to time $N^{\alpha}$, is $o_{P}(1)$ then the limiting
behavior of the rescaled SIS process $X (N^{\alpha}t)/N^{\alpha}$
should be no different from that of the branching envelope $Z
(N^{\alpha}t)/N^{\alpha}$. An easy calculation shows that this will be
the case when $\alpha < 1/2$. When $\alpha =1/2$, the accumulated
attrition over the duration of the branching envelope will be on the
same order of magnitude as the fluctuations, and so the rescaled SIS
process should have a genuinely different asymptotic behavior from the
branching envelope.  Our main result makes this precise:

\begin{theorem}\label{theorem:main}
Assume that process $X (t)=X^{N} (t)$ has infinitesimal transition
probabilities (\ref{eq:rates}). If for some constants $\alpha \leq 1/2$
and $b>0$ the number of individuals initially infected satisfies
$X^{N}(0)\sim bN^{\alpha}$, and if the birth rate (\ref{eq:rates})
satisfies $\beta =\beta_{N}=1+\lambda /N^{\alpha}$, then as $N
\rightarrow \infty$,
\begin{equation}\label{eq:ProcessLimit}
    X^{N} (N^{\alpha}t)/N^{\alpha}
     \stackrel{\mathcal{D}}{\longrightarrow} Y_{t}
\end{equation}
where 
\begin{enumerate}
\item [(a)] if $\alpha <1/2$ then $Y_{t}$ is a Feller diffusion with
drift $\lambda$ and initial state $Y_{0}=b$;
\item [(b)] if $\alpha =1/2$ then $Y_{t}$ is an ``attenuated'' Feller
diffusion with drift $\lambda$ and initial state $Y_{0}=b$, that is,
$Y_{t}$ is a solution to the stochastic differential equation
\begin{equation}\label{eq:attenuatedFellerDiffusion}
    dY_{t}= (\lambda Y_{t}-Y_{t}^{2})\, dt + \sqrt{Y_{t}}dW_{t}
\end{equation}
with $W_{t}$ a standard Wiener process.
\end{enumerate}
\end{theorem}

Note that the ``attenuation'' term $-Y_{t}^{2}$ in the drift of the
limiting process (\ref{eq:attenuatedFellerDiffusion}) can be guessed
from the  form of the attrition $\beta x^{2}/N$. The proof of Theorem
\ref{theorem:main} will be given in section \ref{sec:proof} below.

\section{Size of the Epidemic}\label{sec:size}
The size of an epidemic can be defined as the total number $\xi$ of new
infections during its  entire course. Alternatively, it can be defined
as the total infection time summed over the whole population:
\begin{equation}\label{eq:size}
    S=S_{N}= \int_{0}^{T} X_{t} \,dt.
\end{equation}
Although the two definitions are not the same, it can be shown that
the two quantities have the same asymptotic behavior for large $N$,
that is, $\xi_{N}\sim S_{N}$. (This follows from the fact that the
length of the infection periods for infected individuals are
i.i.d. unit exponential r.v.s.) Because the integral (\ref{eq:size})
is a continuous functional of the path $X_{t}$ (relative to the
Skorohod topology), Theorem \ref{theorem:main} implies that if $X
(0)\sim b\sqrt{N}$ and $\beta =1+\lambda /\sqrt{N}$ then
\begin{equation}\label{eq:sizeLimit}
    S_{N}/N \stackrel{\mathcal{D}}{\longrightarrow}
    \int_{0}^{\tau (0)} Y_{t} \,dt
\end{equation}
where $Y_{t}$ is the attenuated Feller diffusion
(\ref{eq:attenuatedFellerDiffusion}) with initial state $Y_{0}=b$ and
$\tau_{0}$ is the first passage time to $0$ by $Y_{t}$.

By an odd bit of luck, the instantaneous rate $Y_{t}\,dt$ at which
infection time accrues coincides with the rate of change in
accumulated quadratic variation of the semimartingale $Y_{t}$. (Note:
In fact this is really no accident, but rather an artifact of the
fundamental connection between Galton-Watson processes and random
walks via the ``depth-first search'' algorithm. See \cite{aldous} for
more on this.) This suggests making the natural time change to the
diffusion $Y_{t}$ so as to make the instantaneous quadratic variation
constant. The new time scale $s=s (t)$ and the old $t$ are related by
\begin{equation}\label{eq:timeChange}
    ds = Y_{t} \, dt,
\end{equation}
and so $\int Y_{t} \,dt =\int \,ds $ is the limit of the rescaled
epidemic sizes $S_{N}/N$.  The time-changed process
$V_{s}=Y_{t (s)}$ satisfies the stochastic differential equation
\begin{equation}\label{eq:sdeV}
    dV_{s}= (\lambda - V_{s})\,ds + d\tilde{W}_{s}
\end{equation}
where $\tilde{W}_{s}$ is again a standard Wiener process. Setting
$U_{s}=V_{s}-\lambda$, one obtains the stochastic differential
equation for the standard Ornstein-Uhlenbeck process:
\begin{equation}\label{eq:sdeOU}
    dU_{s}=-U_{s}\,ds +d\tilde{W}_{s}
\end{equation}
This proves the following.

\begin{corollary}\label{corollary:size}
If $X (0)\sim b \sqrt{N}$ and $\beta =1+\lambda /\sqrt{N}$ then
\begin{equation}\label{eq:sizeLimitOU}
    S_{N}/N \stackrel{\mathcal{D}}{\longrightarrow}
    \tau (b-\lambda ; -\lambda)
\end{equation}
where $\tau (x;y)$ is the time of first passage to $y$ by a standard
Ornstein-Uhlenbeck process started at $x$.
\end{corollary}

The Laplace transforms of the distributions of $\tau (x;y)$ can be
expressed in terms of parabolic cylinder (Weber) functions: see
Darling and Siegert \cite{darling-siegert}. These do not invert
easily. However, in the special case
$\lambda =0$ (the case corresponding to the critical SIS epidemic!),
the distribution of $\tau (b;0)$ has a simple closed form:
\begin{equation}\label{eq:reflection}
    P\{\tau (b;0)>s \}  = P^{b}\{U_{s}>0 \}-P^{b}\{U_{s}<0 \}.
\end{equation}
This can be obtained from a reflection principle, using the symmetry
of the Ornstein-Uhlenbeck process about the origin.

\section{Proof of Theorem \ref{theorem:main}}\label{sec:proof}

We prove Theorem \ref{theorem:main} using the weak machinery developed
in  Ethier and Kurtz~\cite{et-ku}, which reduces the problem to
checking convergence, in an appropriate sense, of infinitesimal
generators. Write 
\begin{equation}\label{eq:yn}
	Y^{N}_{t}=X^{N} (N^{\alpha}t)/N^{\alpha}
\end{equation}
for the rescaled epidemic process, and denote by $E^{N}_{y}$ the corresponding
expectation operator under the initial condition $Y^{N}_{0}=y$. For
For $f\in\hat C[0,\infty)$, define
\begin{gather}\label{eq:generators}
	\mathcal{G}^Nf(y)=\lim_{h\to0}\frac{E_y^N\big[f(Y_h)-f(y)\big]}{h}
	\quad \text{and}\\
	\mathcal{G}f(y)=(\lambda y)\cdot\frac{\partial f}{\partial
	 y}(y) +y\cdot\frac{\partial^2 f}{\partial y^2}(y) \quad
\text{if} \; \alpha <1/2, \\
	\mathcal{G}f(y)=(\lambda y-y^2)\cdot\frac{\partial f}{\partial
	 y}(y) +y\cdot\frac{\partial^2 f}{\partial y^2}(y) \quad
\text{if}\; \alpha =1/2.
\end{gather}
By Corollary 1.2, sec.~11.2 of \cite{et-ku}, the operator
$\mathcal{G}$ restricted to $\hat{C}[0,\infty)\cap C^{2} (0,\infty)$
generates a Feller semigroup on $\hat C[0,\infty)$, and Proposition
3.3, sec.~1.5 implies that $C_c^\infty[0,\infty)$ is a core for the
generator. (An easy calculation, which we omit, shows that $0$ is an
exit boundary and $\infty$ is a natural boundary in both cases.)
Moreover, the Markov processes determined by these Feller semigroups
can be constructed so as to satisfy the stochastic differential
equations (\ref{eq:feller}) and (\ref{eq:attenuatedFellerDiffusion}),
respectively.  By Theorem 2.5 on page 167 and Theorem 6.1 on page~28,
to prove convergence (\ref{eq:ProcessLimit}) it is enough to show that
for each $f$ in the core of $\mathcal{G}$ the generators converge in
the sense of the following lemma.

\begin{lemma}
\label{lem:genconv} Let $f\in C_c^\infty[0,\infty)$ then
\begin{equation}\label{eq:generatorConvergence}
	 \lim_{N\to\infty}\sup_{y\in[N]\big/\!N^{\alpha}}
	 \left|\mathcal{G}^Nf(y)-\mathcal{G}f(y
	 \right|=0. 
\end{equation}
\end{lemma}

\proof 
Consider first the case $\alpha =1/2$. 
The first step is to calculate $\mathcal{G}^Nf$ for $f\in
C^{\infty}_{c}[0,\infty)$. Using the infinitesimal
transition probabilities (\ref{eq:SIRrates}), we have (with
$x=yN^{1/2}$ and $h=tN^{1/2}$)
\begin{multline*}
	E_y^N\big[f(Y_h^N)-f(y)\big]\\
	 =\left[f\big(1+N^{1/2}\big)-f(y)\right]\times
	\left[\big(1+{\lambda}{{N^{-1/2}}\big)
	 y{N^{1/2}}\big(1-{y}N^{-1/2}}\big)hN^{1/2}\right]\\
	\quad+\left[f\big(1-N^{-1/2}\big)-f(y)\right]\times
	\left[y{N^{1/2}}\cdot hN^{1/2}\right] +o (Nhy)\\
	 =\left[f\big(1+{N^{-1/2}}\big)-f(y)\right]\times
	\left[\big(yN^{1/2}+\lambda y-y^2-
		{\lambda y^2}N^{-1/2}\big)hN^{1/2}\right] \\
	 \quad+\left[f\big(1-N^{-1/2}\big)-f(y)\right]\times
	 \left[yN^{1/2}\cdot h\sqrt{N}\right]+o (Nhy).
\end{multline*}
The error term $o (Nhy)$ is uniform in $y$, because $f$ is assumed to
have compact support. Taking the limit of this expression as $h
\rightarrow 0$ yields
\begin{gather*}
	\mathcal{G}^Nf(y)=
	(f\big(1+N^{-1/2}\big)-f(y))\times N^{1/2}
	\big(\lambda y-y^2-{\lambda y^2}N^{-1/2}\big)\\
	\qquad+\left[ \big(f\big(1-N^{-1/2}\big)-f(y)\big)
	 -\big(f(y)-f(y-N^{-1/2})\big) 
	 \right]\times Ny.
\end{gather*}
Since $f\in C_c^\infty[0,\infty)$, there exists a constant $C>0$
so that $f$ and all its partial derivatives vanish for $y>C$. Therefore,
uniformly in all $y\in[N]\big/\!\sqrt{N}$,
\begin{equation*}
	\lim_{N\to\infty}\mathcal{G}^Nf(y)
	 =\frac{\partial f}{\partial y}(y)\times(\lambda y-y^2)
	 +\frac{\partial^2f}{\partial y^2}(y)\times y =\mathcal{G}f (y).
\end{equation*}

A similar calculation establishes convergence of generators when
$\alpha <1/2$.
\qed

\section{The SIR Epidemic Revisited}\label{sec:sir} The
continuous-time SIR epidemic differs from the SIS epidemic in that
individuals may only be infected once: upon recovery, individuals are
effectively removed from the population. Thus, the state at any time
$t$ is determined by two variables, the number currently infected ($I
(t)=I^{N} (t)$) and the number removed ($R (t)=R^{N} (t)$). These take
valued in the set of nonnegative integer pairs $(i,r)$ such that
$0\leq i+r\leq N$, where $N$ is the (original) population size. The
instantaneous transition rates are as follows:
\begin{align}\label{eq:SIRrates}
(i,r) &\mapsto (i-1,r+1) \quad &\text{at rate} \quad  i\, dt ;\\
\notag 
(i,r) &\mapsto (i+1,r) \quad &\text{at rate} \quad   \beta i (N-i-r)\, dt /N.
\end{align}
All states $(i,r)$ with $i=0$ are absorbing: the epidemic ends at the
first time one of these states is visited.

As for the SIS epidemic, if the numbers of infected and removed
individuals are small compared to the total population size $N$, then
the second transition rate in (\ref{eq:SIRrates}) reduces to $\beta i
\,dt$, and so the process $I (t)$ evolves approxiamtely as the
branching process (\ref{eq:gwe}). Therefore, by the same logic as in
sec.~\ref{ssec:mainResult}, the  limiting behavior of the epidemic can
be deduced by examination of the accumulated attrition over the
duration of the branching process. The result is as follows.

\begin{theorem}\label{theorem:2}
Assume that $(I^{N} (t),R^{N} (t))$ has instantaneous transition rates
(\ref{eq:SIRrates}), and assume that $R^{N} (0)=0$.
If for some $\alpha \leq 1/3$ and $b>0$ the number $I^{N (0)}$ of
individuals initially infected satisfies $I^{N} (0)\sim b N^{\alpha}$,
and if the birth rate $\beta$ satisfies $\beta =1+\lambda
/N^{\alpha}$, then as $N \rightarrow \infty$,
\begin{equation}\label{eq:subT}
	\begin{pmatrix}
	I^{N} (t)\\
	R^{N} (t)
	\end{pmatrix}
\stackrel{\mathcal{D}}{\longrightarrow}
	\begin{pmatrix}
	I (t)\\
	R (t)
	\end{pmatrix}
\end{equation} 
where (i) if $\alpha <1/3$, the limit process $(I (t),R (t))$ satisfies
\begin{align}\label{eq:subTlim}
dI (t)&= \lambda I (t)\,dt +\sqrt{I (t)}\,dW_{t}\\
\notag 
dR (t)&= I (t)\,dt;
\end{align}
and (ii) if $\alpha =1/3$,
\begin{align}\label{eq:Tlim}
dI (t)&= (\lambda I (t) - I (t)R (t))\,dt +\sqrt{I (t)}\,dW_{t}\\
\notag 
dR (t)&= I (t)\,dt.
\end{align}
\end{theorem}

Martin-L\"{o}f's result (\ref{eq:ml}) can be easily recovered from
Theorem \ref{theorem:2} by the same device as used in
sec.~\ref{sec:size} above. Define the new time scale $s$ by $ds =
I_{t}\,dt$ and corresponding time-changed process $dJ (s)= dI
(t)$. Then the total size of the epidemic is just the integral $\int
ds$ up to the time of first passage to $0$ by $J (s)$. But $J (s)$ is
just the Wiener process with a quadratic drift, and so (\ref{eq:ml}) follows.

Theorem \ref{theorem:2} can be proved either by martingale methods or
by use of the Ethier-Kurtz machinery. The latter approach is mildly
complicated by the fact that the generator
\begin{equation}\label{eq:IRGenerator}
	\mathcal{G}= (\lambda i -ir)\partial_{i}+\sqrt{i}\partial_{ii}
		+i\partial_{r}
\end{equation}
is not elliptic, but rather parabolic, and singular along the $i=0$
axis. The singularity at $i=0$ can be handled by truncating the state
space: To prove (\ref{eq:subT}), it suffices to prove weak convergence
for the processes $(I^{N} (t\wedge \tau_{\varepsilon }),R^{N} (t\wedge
\tau_{\varepsilon}))$, where $\tau_{\varepsilon }$ is the time of
first passage to the level $i=\varepsilon$. Nonellipticity of the
generator may be handled by using standard existence results from the
theory of parabolic PDE (see \cite{friedman}) to verify the hypotheses
of the Hille-Yosida theorem (Th.~2.2 of \cite{et-ku}). Weak
convergence of the truncated processes may then be proved by checking
convergence of generators; this is another routine calculation similar
to that carried out for the SIS epidemic in sec.~\ref{sec:proof}
above.

\end{document}